\documentclass[10 pt, conference]{ieeeconf}

\usepackage{graphicx}
\usepackage{dsfont} 
\usepackage{amsmath}
\usepackage{amssymb}
\usepackage{url}

\usepackage{color}
\usepackage{subcaption}

\newcommand{\mc}{\mathcal}
\newcommand{\doi}[1]{\href{http://dx.doi.org/#1}{\normalsize{\textsc{doi:}}~\nolinkurl{#1}}}
\newcommand{\arxiv}[1]{\href{http://arxiv.org/abs/#1}{\normalsize{\textsc{arxiv:}}~\nolinkurl{#1}}}

\renewcommand{\epsilon}{\varepsilon}
\renewcommand{\phi}{\varphi}

\newcommand{\R}{\mathbb{R}}

\let\originalleft\left
\let\originalright\right
\renewcommand{\left}{\mathopen{}\mathclose\bgroup\originalleft}
\renewcommand{\right}{\aftergroup\egroup\originalright}

\def\clap#1{\hbox to 0pt{\hss#1\hss}}

\newcommand{\norm}[1]{\left\lVert #1\right\rVert}

\newcommand{\set}[1]{\left\{ #1\right\}}

\let\p=\paren

\newcommand{\sqparen}[1]{\left[ #1\right]}

\let\sp=\sqparen

\newcommand{\derivoper}{\mathrm{d}}

\newcommand{\deriv}[3][]{%
  \ifx\relax#1\relax{
    \frac{\derivoper #2}{\derivoper #3}%
  }\else{%
    \frac{\derivoper^{#1} #2}{\derivoper #3^{#1}}%
  }\fi%
}
\newcommand{\pderiv}[3][]{%
  \ifx\relax#1\relax{
    \frac{\partial #2}{\partial #3}%
  }\else{%
    \frac{\partial^{#1} #2}{\partial #3^{#1}}%
  }\fi%
}
\newcommand{\grad}[2][]{%
  \ifx\relax#1\relax{
    \nabla #2%
  }\else{%
    \nabla_{#1} #2%
  }\fi%
}
\newcommand{\diver}[2][]{%
  \ifx\relax#1\relax{
    \nabla \cdot #2%
  }\else{%
    \nabla_{#1} \cdot #2%
  }\fi%
}
\newcommand{\curl}[2][]{%
  \ifx\relax#1\relax{
    \nabla \times #2%
  }\else{%
    \nabla_{#1} \times #2%
  }\fi%
}
\newcommand{\lapl}[2][]{%
  \ifx\relax#1\relax{
    \Delta #2%
  }\else{%
    \Delta_{#1} #2%
  }\fi%
}

\newcommand{\smats}[1]{%
  \sqparen{\begin{smallmatrix}#1 \end{smallmatrix}}%
}

\IEEEoverridecommandlockouts

\usepackage{cleveref}
\newcommand{\hs}{\mathcal{H}}

\newcommand{\M}{\mathcal{M}}
\newcommand{\Me}{\mathcal{M}^\epsilon}
\newcommand{\J}{\mathcal{J}}
\newcommand{\BVU}{PC\p{\sp{0,T}, U}}
\newtheorem{assumption}{Assumption}
\newtheorem{definition}{Definition}
\newtheorem{lemma}{Lemma}
\newtheorem{theorem}{Theorem}
\newtheorem{proposition}{Proposition}
\usepackage{comment}

\usepackage{mathrsfs}

\DeclareMathAlphabet{\mathpzc}{OT1}{pzc}{m}{it}

\begin{document}
\title{%
 A New Solution Concept and\\  Family of Relaxations for Hybrid Dynamical Systems}
\author{Tyler Westenbroek \and Humberto Gonzalez \and S. Shankar Sastry
\thanks{T. Westenbroek  and S. S. Sastry are with the Department of Electrical Engineering and Computer Sciences, University of California, Berkeley.  Emails: {\tt\footnotesize $\{$westenbroekt, sastry$\}$@eecs.berkeley.edu}. Humberto Gonzalez is with BlackThorn Therapeutics. Email: hgonzale@gmail.com. This work was supported by FORCES (Foundations Of Resilient CybEr-physical Systems), National Science Foundation award number CNS-1239166, and HICON-LEARN (design of HIgh CONfidence LEARNing-enabled systems), Defense Advanced Research Projects Agency award number FA8750-18-C-0101.
}
}

\maketitle

\begin{abstract}
We introduce a holistic framework for the analysis, approximation and control of the trajectories of hybrid dynamical systems which display event-triggered discrete jumps in the continuous state. We begin by demonstrating how to explicitly represent the dynamics of this class of systems using a single piecewise-smooth vector field defined on a manifold, and then employ Filippov's solution concept to describe the trajectories of the system. The resulting \emph{hybrid Filippov solutions} greatly simplify the mathematical description of hybrid executions, providing a unifying solution concept with which to work. Extending previous efforts to regularize piecewise-smooth vector fields, we then introduce a parameterized family of smooth control systems whose trajectories are used to approximate the hybrid Filippov solution numerically. The two solution concepts are shown to agree in the limit, under mild regularity conditions.
\end{abstract}


\section{Introduction}

Hybrid dynamical systems are a natural abstraction for many physical and cyber-physical systems \cite{hiskens2000trajectory}, \cite{grizzle2001asymptotically}. Yet, despite extensive efforts to characterize the subtle interactions that arise between continuous and discrete dynamics \cite{goebel2004hybrid}, \cite{ames2005homology}, \cite{lygeros2003dynamical}, \cite{hiskens2000trajectory}, the dynamical properties of hybrid systems are not adequately understood. 

In this paper, we extend the geometric approach advocated in \cite{simic2005towards}, and reduce a given hybrid system to a piecewise-smooth control system defined on a manifold. Applying Filippov's convention \cite{filippov2013differential} to this system, we characterize the trajectories of each hybrid system using a single differential inclusion. These \emph{hybrid Filippov solutions} are not defined using discrete transitions, and thus are not subject to many of the theoretical challenges facing traditional constructions of hybrid executions. To illustrate this point, throughout the paper we will consider how our approach simplifies the mathematical description of Zeno executions \cite{johansson1999regularization}, hybrid trajectories which undergo an infinite number of discrete transitions in a finite amount of time. 

However, numerical tools for simulating, analyzing and controlling Filippov systems are far less developed than they are for smooth dynamical systems. Therefore, we build on the work of \cite{llibre2008sliding} and demonstrate how to approximate the hybrid Filippov solution using a parameterized family of smooth, stiff control systems defined on the relaxed topology from \cite{burden2015metrization}. Under standard regularity assumptions, these relaxations are shown to recover the hybrid Filippov solution as the appropriate limit is taken. Our relaxations provide a foundation for controlling hybrid systems using standard numerical algorithms developed for smooth control systems \cite{polak2012optimization}, \cite{schattler2012geometric}. Proofs of claims made in the paper can be found in \cite{arxiv}, where numerous examples are explored and a thorough literature review is provided.



\section{Mathematical Notation}\label{sec:math}

In this section we fix mathematical notation used throughout the paper. Given a set $D$, $\partial D$ is the boundary of $D$ and $int(D)$ is the interior of $D$. For a topological space $V$, we let $\mc{B}(V)$ denote all subsets of $V$. Given a metric space $(X, d)$, we denote the ball of radius $\delta$ centered at $x \in X$ by $B^{\delta}(x)$. The 2-norm is our metric of choice for finite-dimensional real spaces, unless otherwise noted. We use $\overline{co}S$ to denote the convex closure of a set $S$, which is a subset of some vector space $V$. The \emph{disjoint union} of a collection of sets $\set{D_j}_{j \in \mc{J}}$  is denoted ${\coprod_{j \in \mc{J}} D_j = \bigcup_{j \in \mc{J} }D_j \times \set{j}}$, which is endowed with the piecewise topology. Given $x \in D_{j} \times \set{j}$, we will frequently abuse notation and simply write $x \in D_j$ when context makes our meaning clear. Throughout the paper we use the term \emph{smooth} to mean infinitely differentiable and it is understood that \emph{diffeomorphisms} are smooth mappings. 

We assume familiarity with the notions of \emph{topological manifolds} and \emph{quotient spaces}, and refer the reader to \cite{arxiv} or \cite{manifolds2012introduction} if they are unfamiliar with these concepts. For a topological space $\mathcal{S}$ and a function $f \colon A \to B$,  where $A,B \subset \mathcal{S}$, we define $\Lambda_f$ to be the smallest equivalence relation containing  the set $\set{(a,b) \in \mathcal{S}\times \mathcal{S} \colon a \in f^{-1}(b)}$, and denote the set of equivalence classes by $ \frac{\mathcal{S}}{\Lambda_f}$. There is a natural \emph{quotient projection} $\pi \colon \mc{S} \to \frac{S}{\Lambda_f}$ taking each $s \in S$ to its equivalence class $[s] \in \frac{S}{\Lambda_f}$ and we endow $ \frac{\mathcal{S}}{\Lambda_f}$ with the  quotient topology. Given a smooth manifold $\M$ and compact set of allowable inputs $U \subset \R^m$, we say that a function $F \colon \M \times U \to T \M$ is a vector field on $\M$ if, for each $X \in \M$ and $u \in U$, $F(X,u) \in T_{X}\M$. Throughout the paper, we will consider input signals in the space of piecewise-continuous controls, which will be denoted $\BVU$. Given a vector field defined on $\R^n$, we will frequently abuse notation and say $f \colon \R^n \times U \to \R^n$ when we should write $f \colon \R^n \times U \to T\R^n$.


\section{Filippov Solutions}\label{sec:filippov}

We now introduce Filippov's solution concept \cite{filippov2013differential} for differential equations with discontinuous right-hand sides. In the sequel we will use these solutions to locally describe the dynamics of our class of hybrid systems. Throughout the paper, we will primarily focus on the existence and uniqueness of solutions. However, given our construction of the hybrid Filippov solution, it is straightforward to apply any known property of Filippov's solution concept to our class of hybrid systems. 

To simplify exposition, throughout the majority of the paper we will restrict our attention to hybrid trajectories that can be locally described using a bimodal discontinuous vector field $f\colon \R^n \times U \to \R^n$ of the form
\begin{equation}\label{eq:pws_system}
f(x,u) = \begin{cases}
f_1(x,u) & \text{ if } x \in D_1\\
f_2(x,u) & \text{ if } x \in D_2,
\end{cases}
\end{equation}
where, given a smooth, regular function $g \colon \R^n \to \R$, we put $D_1 = \set{x \in \R^n \colon g(x) <0}$ and $D_2 = \set{x \in \R^n \colon g(x)>0}$, and for $i \in \set{1,2}$ we require $f_i \colon \R^n \times U \to \R^n$ to be smooth and globally Lipschitz continuous. Note that the discontinuity set for $f$ is confined to the surface ${\Sigma \colon = \set{x \in \R^n \colon g(x)= 0}}$. In Section \ref{sec:examples} and in \cite{arxiv} we consider several examples of hybrid models whose dynamics can be locally represented by piecewise-smooth vector fields with multiple, overlapping surfaces of discontinuity. 

The \emph{Filippov Regularization} of a general discontinuous vector field $\tilde{f} \colon \R^n \times U \to \R^n$ is the set-valued map $\mc{F}[\tilde{f}] \colon \R^n \times U \to \mc{B}(\R^n)$ defined by
\begin{equation}\label{eq:filippov}
\mc{F}\sp{\tilde{f}}(x, u)  = \bigcap_{\delta > 0} \bigcap_{\mu (S) = 0}  \overline{co} f( B^\delta(x) \setminus S, u),
\end{equation}
where $\bigcap_{\mu (S) = 0}$ denotes the intersection over all sets of zero measure in the sense of Lebesgue. We say that a \emph{Filippov solution} for the differential equation $\dot{x}(t) = \tilde{f}(x(t),u(t))$  given initial data $x_0 \in D$ and control $u \in \BVU$ is an absolutely continuous curve $x \colon \sp{0,T} \to \R^n$ such that $x(0) = x_0$ and $\dot{x}(t) \in \mc{F} \sp{\tilde{f}}(x(t), u(t)) \text{ a.e. } t \in \sp{0,T}$.

We next consider the existence and uniqueness of Filippov solutions for the piecewise-smooth vector field \eqref{eq:pws_system}. Lemmas \ref{lemma:existence} and \ref{lemma:uniqueness} are adapted from results in \cite[Chapter 2]{filippov2013differential}, and we discuss how to obtain these particular results in \cite{arxiv}.
\begin{lemma}\label{lemma:existence}
Consider the discontinuous system \eqref{eq:pws_system}. For each $(x_0, u)  \in  \R^n \times \BVU$ there exists a Filippov solution $x \colon \sp{0,T} \to \R^n$ for the differential equation ${\dot{x}(t) = f(x(t),u(t))}$.
\end{lemma}

That is, Filippov solutions for the discontinuous system \eqref{eq:pws_system} exist on bounded time intervals. The following is a sufficient condition for the uniqueness of Filippov solutions corresponding to \eqref{eq:pws_system}.
\begin{assumption}\label{ass:transverse}
Consider the discontinuous system \eqref{eq:pws_system}. For each $(x,u) \in \Sigma \times U$ either ${\nabla g(x) \cdot f_1(x,u) > 0}$ or ${\nabla g(x) \cdot f_2(x,u) < 0}$.  
\end{assumption}

Assumption \ref{ass:transverse} rules out a number of pathological cases where trajectories skim the surface of discontinuity at points where $f_1(x,u)$ and $f_2(x,u)$ are both tangent to $\Sigma$.

\begin{lemma}\label{lemma:uniqueness}
Let Assumption \ref{ass:transverse} hold for the discontinuous system \eqref{eq:pws_system}. Then for each $(x_0,u)  \in \R^n \times \BVU$ there is a unique Filippov solution $x \colon \sp{0,T} \to \R^n$ for the differential equation $\dot{x}(t) = f(x(t),u(t))$.
\end{lemma}

In cases where $x \in \Sigma$, ${\nabla g(x) \cdot f_1(x,u) >0}$ and ${\nabla g(x) \cdot f_2(x,u) <0}$, Filippov's convention allows us to equivalently describe the dynamics of the system using a \emph{sliding vector field} defined on $\Sigma$ \cite[Chapter 1]{filippov2013differential}.
\section{Relaxed Filippov Systems}\label{sec:smooth_systems}
Teixeira has developed a framework for approximating the dynamics of autonomous, piecewise-smooth vector fields using a parameterized family of smooth, stiff vector fields (see e.g. \cite{llibre2008sliding}).  Here, we generalize the approach to control systems of the form \eqref{eq:pws_system}. In Section \ref{sec:relaxed_hybrid_sys}, we will modify this approach to locally approximate the hybrid Filippov solution. The following class of functions will be used to regularize \eqref{eq:pws_system} near the surface of discontinuity. 

\begin{definition}
We say that $\phi \in C^{\infty}\p{\R, \sp{0,1}}$ is a \emph{transition function} if $i)$ $\phi(a) = 0$ if $a \leq -1$, $ii)$ $\phi(a) = 1$ if $a \geq 1$, and $iii)$ $\phi$ is monotonically increasing on $\p{-1, 1}$.
\end{definition}

For the rest of the section, we assume that a single transition function $\phi$ has been chosen. We then define for each $\epsilon >0$ the function $\phi^\epsilon \colon \R^n \to \R$ by $\phi^\epsilon(x) = \phi\p{\frac{g(x)}{\epsilon}}$, and then define the $\epsilon$-relaxation of \eqref{eq:pws_system} to be ${f^\epsilon \colon \R^n \times U \to \R^n}$ where
\begin{equation}\label{eq:smooth_filippov}
f^\epsilon(x,u) = \p{1-\phi^\epsilon(x)} f_1(x,u) + \phi^\epsilon(x)f_2(x,u).
\end{equation}
The relaxation occurs by smoothing \eqref{eq:pws_system} along
\begin{equation}
\Sigma^\epsilon\colon = \set{x \in \R^n \colon - \epsilon \leq g(x) \leq \epsilon}.
\end{equation}
Indeed, note that $f^\epsilon(x,u) = f_1(x,u)$ if $x \in D_1 \setminus \Sigma^\epsilon$, $f^\epsilon(x,u) = f_2(x,u)$ if $x \in D_2 \setminus \Sigma^\epsilon$, and $f^\epsilon(x,u)$ produces a convex combination of $f_1(x,u)$ and $f_2(x,u)$ if $x \in \Sigma^\epsilon$. 

\begin{lemma}
For each $\epsilon >0$, the relaxed vector field \eqref{eq:smooth_filippov} is smooth. 
\end{lemma}

The following two theorems characterize the behavior of the relaxed system as $\epsilon \to 0$. In Section \ref{sec:relaxed_hybrid_sys}, these results will be extended to our relaxed solution concept for hybrid systems.
\begin{theorem}\label{thm:convergence1}
Let Assumption \ref{ass:transverse} hold for the discontinuous system \eqref{eq:pws_system}. Fix $(x_0, u ) \in \R^n \times \BVU$ and let $x \colon \sp{0,T} \to \R^n$ be the corresponding (unique) Filippov solution for \eqref{eq:pws_system}, and for each $\epsilon>0$ let $x^\epsilon \colon \sp{0,T} \to \R^n$ be the corresponding solution to the relaxed system \eqref{eq:smooth_filippov}. Then $\exists C >0$ and $\epsilon_0 >0$ such that for each $\epsilon \leq \epsilon_0$
\begin{equation}
\norm{x - x^\epsilon}_{\infty} \leq C \epsilon.
\end{equation}
\end{theorem}

The result is obtained by transforming the relaxed system into a \emph{singular perturbation problem} (see e.g. \cite[Chapter 6.3]{sastry2013nonlinear}), and applying standard convergence results from the literature. However, even when Assumption \ref{ass:transverse} is not satisfied, the solutions of the relaxed system still converge uniformly to a well-defined limit. 

\begin{theorem}\label{thm:convergence2}
Fix $(x_0,u) \in \R^n \times \BVU$ and for each $\epsilon >0$ let $x^\epsilon \colon \sp{0,T} \to \R^n$ be the corresponding solution to the relaxed system \eqref{eq:smooth_filippov}. Then there exists an absolutely continuous curve $x^0 \colon \sp{0,T} \to \R^n$ such that
\begin{equation} \vspace{10pt}
\lim_{\epsilon \to 0} \  \norm{x^\epsilon - x^0}_{\infty} = 0.
\vspace{10pt}
\end{equation}
\end{theorem}



\section{Hybrid Dynamical Systems}\label{sec:hybrid_sys}
In this section we define our class of hybrid dynamical systems, discuss technical assumptions made throughout the paper, and introduce and characterize the topological spaces upon which we will define the hybrid Filippov solution and its relaxations. 

\subsection{Hybrid Dynamical Systems}
\begin{definition} \label{def:hybrid_sys}
A hybrid dynamical system is a seven-tuple 
\begin{equation}
\mathcal{H} =  \p{\mathcal{J}, \Gamma, \mathpzc{D}, U, \mathpzc{F}, \mathpzc{G}, \mathpzc{R}},
\end{equation}where:
\begin{itemize}
\item $\mathcal{J}$ is a finite set indexing the discrete states of $\mathcal{H}$; 
\item $\Gamma \subset \mathcal{J} \times \mathcal{J}$ is the set of edges, forming a graphical structure over $\mathcal{J}$, where edge $e = (j,j') \in \Gamma$ corresponds to a transition from $j$ to $j'$;
\item $\mathpzc{D} =  \set{D_j}_{j \in \mathcal{J}}$ is the set of domains, where $D_j \subset \R^n$ is a smooth, connected $n$-manifold with boundary;
\item $U \subset \R^m$ is a compact set of allowable inputs;
\item $\mathpzc{F} = \set{f_j}_{j \in \mathcal{J}}$ is the set of vector fields, where each $f_j \colon \R^n \times U \to \R^n$ defines the dynamics on $D_j$; 
\item $\mathpzc{G} = \set{G_e}_{e \in \Gamma}$ is the set  of guards, where each $G_{(j,j')} \subset \partial D_j$ is a smooth, embedded $(n-1)$-manifold; 
\item $\mathpzc{R} = \set{R_e}_{e \in \Gamma}$ is the set of reset maps, where $R_{(j,j')} \colon \R^n \to \R^n$ and $R_{(j,j')}(G_{(j,j')}) \subset \partial D_{j'}$. 
\end{itemize}
\end{definition} 

Before enumerating the technical assumptions we make throughout the paper, we provide a simple definition for the executions of a hybrid dynamical system that resembles most definitions in the literature, which we will use to highlight the features of the hybrid Filippov solution and its relaxations.

\begin{definition}\label{def:classical_hybrid_traj}
Let $\hs$ be a hybrid dynamical system and let $(x_0,u) \in D_{q_1} \times \BVU$. Let $\set{q_{k}}_{k=1}^N$ be a sequence of discrete states, where $N \in \mathbb{N}$ is possibly infinite, let $\set{e_{k}}_{k=1}^{N-1}$ be a sequence of edges where $e_k = (q_k, q_{k+1})$, let $0 =t_1 \leq t_2 \leq \dots \leq t_{N} \leq t_{N+1} =T$, and let $\set{x_k}_{k=1}^{N}$ where $x_k \colon \sp{t_k, t_{k+1}} \to D_{q_k}$ is absolutely continuous satisfy the following conditions: 
\begin{align}\label{eq:fbvp1}
x_1(0) =& x_0\\
 \dot{x}_k(t) = f_{q_k}(x_k(t)&,u(t)), \ \forall t \in [t_k, t_{k+1}) \\
x_{k}(t_{k+1}) & \in G_{e_k} \\
x_{k+1}(t_{k+1}) = R_{e_k}&(x_{k}(t_{k+1})).
\end{align}
Then we say that $x \colon \sp{0,T} \to \coprod_{j \in \J} D_j$ where
\begin{equation}
x(t) = x_k(t) \times \set{q_{k}}, \ \  \forall t \in \sp{t_k,t_{k+1}}
\end{equation}
is a \emph{hybrid execution} corresponding to $(x_0,u)$.
\end{definition}

\begin{figure}[t!]
  \centering
  \vspace{10pt}
  \includegraphics[width= .5\columnwidth,  keepaspectratio ]{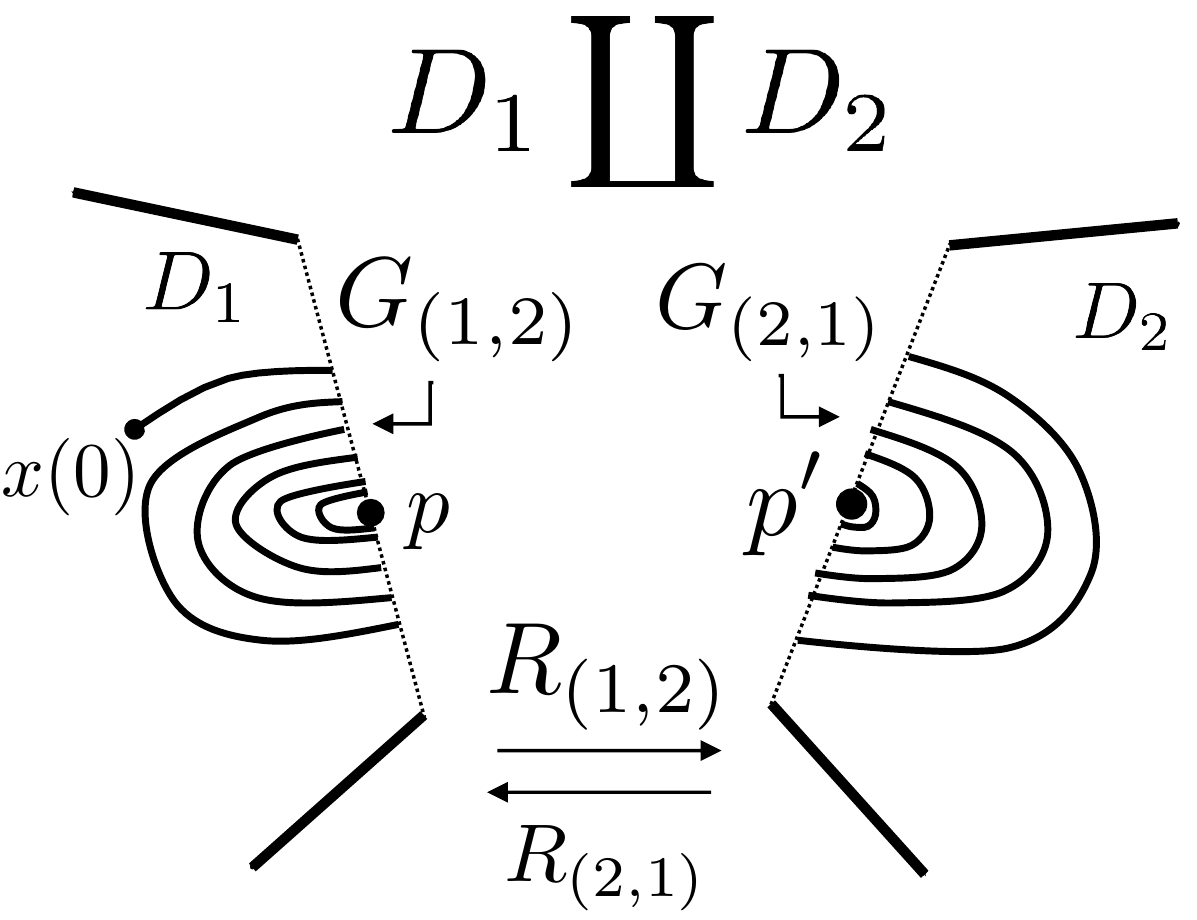}%
  \caption{A hybrid execution $x$ evolves from initial condition $x(0)$ on the disjoint union $D_1 \coprod D_2$. The trajectory is Zeno, undergoing an infinite number of discrete jumps in a finite amount of time. The portions of $x$ in $D_1$ accumulate to the point $p$, and the portions in $D_2$ accumulate to the point $p'$. The hybrid system has two edges: $(1,2)$ and $(2,1)$, where $R_{(1,2)}(G_{(1,2)})= G_{(2,1)}$ and $R_{(2,1)} = R_{(1,2)}^{-1}$. }
  \label{fig:classic_traj}
  \vspace{-10pt}
  \end{figure}

Note that we have defined a hybrid execution to be multivalued at the transition times $\set{t_{k}}_{k=2}^{N}$. We say that a hybrid execution is \emph{Zeno} if $N = \infty$.  A hybrid execution which is Zeno is depicted in Figure \ref{fig:classic_traj}. The following assumption guarantees the existence and uniqueness of solutions to the ordinary differential equations defined on each continuous domain.
\begin{assumption}\label{ass:lipschitz_vectorfield}
For each $j \in \J$, the vector field $f_j$ is smooth and Lipschitz continuous. 
\end{assumption}

The hybrid Filippov solution will be defined on the \emph{hybrid quotient space} or \emph{hybridfold} from \cite{simic2005towards}, and its relaxations will be defined on the relaxed version of this topology from \cite{burden2015metrization}. We introduce these important concepts in Sections \ref{sec:hybrid_quotient_space} and \ref{sec:relaxed_topology}, but first make assumptions which will ensure that these spaces are sufficiently regular, and also simplify the initial introduction of our framework. We first make a strong assumption about the geometry of guard sets and their images under reset maps.
\begin{assumption}\label{ass:affine}
For each $e =(j,j') \in \Gamma$, $G_e$ and $R_e(G_e)$ are subsets of $(n-1)$-dimensional hyperplanes. Specifically, there exists unit vectors $\hat{g}_e, \hat{r}_e \in \R^n$ and scalars $c_e, d_e$ such that $G_e \subset \tilde{G}_e$ and $R_e(G_e) \subset R_e(\tilde{G}_e) \subset \tilde{R}_e$, where
\begin{enumerate}
\item $\tilde{G}_e  \colon = \set{x \in \R^n \colon g_e(x) \colon = \hat{g}_e^Tx -c_e = 0}$, \text{ and }
\item $\tilde{R}_e  \colon = \set{x \in \R^n \colon r_e(x) \colon = \hat{r}_e^Tx -d_e = 0}$.
\end{enumerate} 
Furthermore, $g_e(x) < 0$ for each $x \in D_{j} \setminus G_e$, and $r_e(x)  > 0$ for each $x \in D_{j'} \setminus R_e(G_e)$.
\end{assumption}

In Sections \ref{sec:coordinate_charts} and \ref{sec:relaxed_topology} we will demonstrate how to explicitly construct a collection of coordinate charts for the hybrid quotient space and its relaxations, and our approach will rely on this assumption. In practice, we have found that it is often possible to get around this assumption by adding auxiliary continuous states to the hybrid system, or by simply choosing a coordinate system in which the assumption is satisfied. \footnote{Since each continuous domain is assumed to be a $n$-dimensional manifold with boundary in Definition \ref{def:hybrid_sys}, for each $e \in \Gamma$ there must exist a collection of boundary charts covering $G_e$ and $R_e(G_e)$. Moreover, in coordinates both of these sets will be defined locally by $\set{(x_1, \dots , x_n) \in \R^n \colon x_n =0}$. Thus, in principle we could always satisfy Assumption \ref{ass:affine} by working in local coordinates where appropriate throughout the paper. However, we make Assumption \ref{ass:affine} since in practice we may not have access to the necessary charts. }  This point is illustrated by our examples in \cite{arxiv}.

The following assumption is crucial for ensuring that the hybrid quotient space and its relaxations are smooth manifolds (see \cite[Theorem 1]{simic2005towards} or \cite[Theorem 3]{burden2015model}). 
\begin{assumption}\label{ass:diff_rmaps}
For each $e  \in \Gamma$ the map $R_e$ is a diffeomorphism. Furthermore, $\nabla R_e$ and $\nabla R_e^{-1}$ are both globally Lipschitz continuous. 
\end{assumption}

Requiring the invertibility of each reset map is a strong technical assumption. A number of application domains, such as robotic bipedal walking \cite{grizzle2001asymptotically}, typically utilize hybrid models with rank-deficient reset maps. Nevertheless, as we illustrate with our examples, it is often possible to transform these hybrid systems into ones which do satisfy Assumption \ref{ass:diff_rmaps} by adding auxiliary continuous states to the system. Oftentimes, these extra states have physically meaningful interpretations. The assumption on the gradients of the reset maps will ensure that the vector fields we define on the (relaxed) hybrid quotient space are sufficiently regular for our purposes. 

Next, we impose an assumption which will ensure that we can locally describe the trajectories of our hybrid systems using vector fields of the form \eqref{eq:pws_system}.
\begin{assumption}\label{ass:overlapping_guards}
The elements of ${\set{G_e}_{e\in \Gamma} \cup \set{R_e(G_e)}_{e \in \Gamma}}$ are mutually disjoint.
\end{assumption}

We initially assume elements of ${\set{G_e}_{e\in \Gamma} \cup \set{R_e(G_e)}_{e \in \Gamma}}$ do not intersect for notational convenience, and will indicate throughout the paper and in \cite{arxiv} ways to weaken this assumption.  Just as in \cite[Theorem 3]{burden2015model}, our final assumption ensures the hybrid quotient space and its relaxations are topological manifolds \emph{without boundary}.
\begin{assumption}\label{ass:boundary}
For each $j \in \J$ and $x \in \partial D_j$ there exists $e \in \Gamma$ such that either $x \in G_e$ or $x \in R_e(G_e)$.
\end{assumption}

Assumption \ref{ass:boundary} is also made primarily for  convenience, as working with manifolds with boundary requires additional overhead. As discussed in \cite{arxiv}, the results of this paper go through in a natural way when this assumption is lifted.

\subsection{The Hybrid Quotient Space}\label{sec:hybrid_quotient_space}

\begin{figure}[t!]
\vspace{10pt}
  \centering
  \includegraphics[width= .7\columnwidth,  keepaspectratio ]{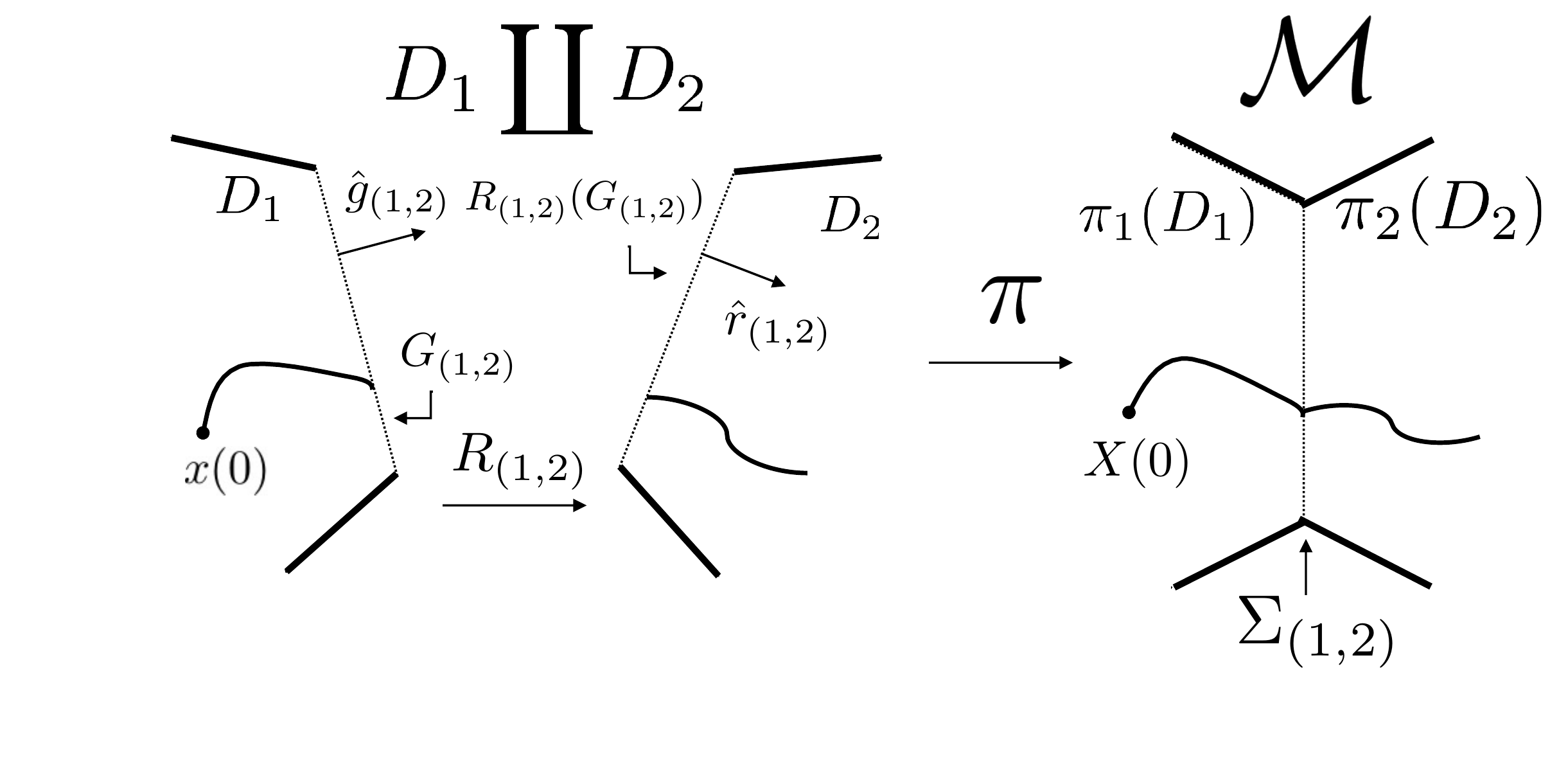}%
  \caption{Construction of the hybrid quotient space from the disjoint union of the continuous domains for a bimodal hybrid system with a single edge $e= (1,2)$. A hybrid execution $x$ transitions from mode 1 to mode 2 on $D_1 \coprod D_2$. The continuous curve $X = \pi \circ x$ is the representation of $x$ on $\M$. }
  \label{fig:hybrid_quotient_space}
    \vspace{-15pt}
\end{figure}

The main idea behind the construction of the hybrid quotient space, which is depicted in Figure \ref{fig:hybrid_quotient_space}, is to identify or "glue" each point $x \in G_e$ to the point $R_e(x) \subset R_e(G_e)$. This process unifies the domains $\set{D_j}_{j \in \J}$ into a single topological space. As depicted in Figure \ref{fig:hybrid_quotient_space}, we can interpret hybrid executions as continuous curves on the hybrid quotient space. We refer the reader to \cite{simic2005towards} for a number of simple examples which clearly illustrate these concepts. In what follows, we closely follow the notation developed in \cite{burden2015metrization}. Formally, for a given hybrid system $\hs$, we define ${\hat{R} \colon \coprod_{e \in \Gamma} G_e \to \coprod_{j  \in \J}  D_j}$ by $\hat{R}(x) = R_e(x)$ if $x \in G_e$, and then define the \emph{hybrid quotient space} to be 
\begin{equation}
\M = \frac{\coprod_{j \in \J} D_j}{\Lambda_{\hat{R}}}.
\end{equation}

Recalling our notation from Section \ref{sec:math}, $\Lambda_{\hat{R}}$ is an equivalence relation on $\coprod_{j \in \J} D_j$, wherein the equivalence class for each point $x \in G_e$ is the set $\set{x , R_e(x)}$. The hybrid quotient space is endowed with a quotient map $\pi \colon \coprod_{j \in \J} D_j \to \M$ that takes each point $\hat{x} \in \coprod_{j \in \J} D_j$ to its equivalence class $[\hat{x}] \in \M$. Since the points $x$ and $R_e(x)$ above belong to the same equivalence class, they are sent to a single point $\pi(x) = \pi(R_e(x)) \in \M$. 

For notational clarity later on, for each $j \in \J$ we define the map $\pi_j \colon D_j \to \M$ by $\pi_j(x) = \pi(x \times \set{j})$, which takes each point in $D_j \subset \R^n$ to the corresponding point in $\M$. As shown in Figure \ref{fig:hybrid_quotient_space}, each domain $D_j$ is represented on $\M$ by the regular domain $\pi_j(D_j) \subset \M$. For each $e =(j,j') \in \Gamma$, $G_e \subset \partial D_j$ and $R_e(G_e) \subset \partial D_{j'}$ are collapsed to
\begin{equation}
\Sigma_{e} \colon = \pi_j(G_e) = \pi_{j'}(R_e(G_e)),
\end{equation}
which is a co-dimension-1 sub-manifold of $\M$ separating the interiors of $\pi_{j}(D_j)$ and $\pi_{j'}(D_{j'})$. \footnote{Consider the hybrid system depicted in Figure \ref{fig:classic_traj}. When constructing the hybrid quotient space for this system, the surfaces $G_{(1,2)}$ and $G_{(2,1)}$ are sent to a single surface. Thus, we can remove one of the edges from this hybrid system to satisfy Assumptions \ref{ass:affine}, without changing the structure of $\M$. As noted in Section \ref{sec:hybrid_filippov}, removing this redundant edge will not affect our two solution concepts, which still capture the behavior of this system. }

As noted in \cite{simic2005towards}, $\M$ is always metrizable when our standing assumptions are satisfied; that is, there exists a state-space metric $d \colon \M \times \M \to \R_+$. One such metric is explicitly constructed in \cite{burden2015metrization}. For our purposes, we assume that a specific metric has been chosen, which we will refer to as $d$. Whenever we refer to a curve as being (absolutely) continuous on $\M$, it is understood that we mean (absolutely) continuous with respect to $d$. Using the map $\pi$, we can descend a hybrid execution to a continuous curve on $\M$, as depicted in Figure \ref{fig:hybrid_quotient_space}. The hybrid Filippov solution will directly generate continuous curves on $\M$.

\subsection{Charting the Hybrid Quotient Space}\label{sec:coordinate_charts}

Under our standing assumptions, the hybrid quotient space is in fact a smooth topological manifold \cite[Theorem 3]{burden2015model}. By regarding $\M$ as a manifold, in Section \ref{sec:hybrid_filippov} we will be able to define a single piecewise-smooth vector field on $\M$ that captures the dynamics of the hybrid system. In this section, we demonstrate how to explicitly construct a set of smoothly compatible coordinate charts for $\M$. 

\begin{figure}[t!]
\vspace{10pt}
  \centering
  \includegraphics[width= .9\columnwidth, height = .95\textheight, keepaspectratio ]{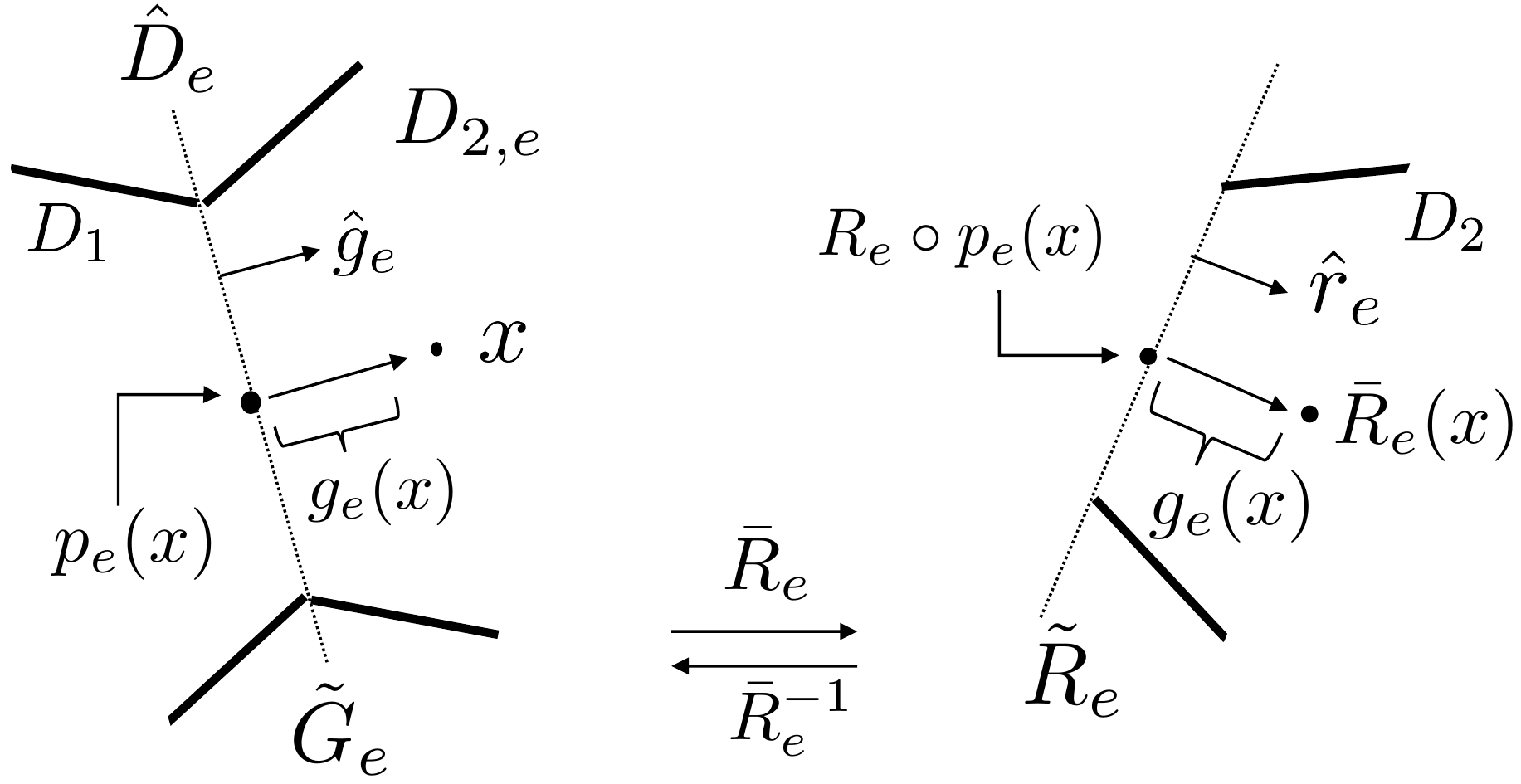}%
  \caption{The domain $D_2$ is smoothly attached to domain $D_1$, using the map $\bar{R}_e$ and resulting in $D_{2,e} \colon = \bar{R}_e^{-1}(D_2)$, for edge $e=(1,2)$. The various components of $\bar{R}_e$ are illustrated. }
  \label{fig:attach}
  \vspace{-20pt}
\end{figure}

\begin{figure}[t]
\vspace{20pt}
  \centering
  \includegraphics[width= .9\columnwidth, height = .95\textheight, keepaspectratio ]{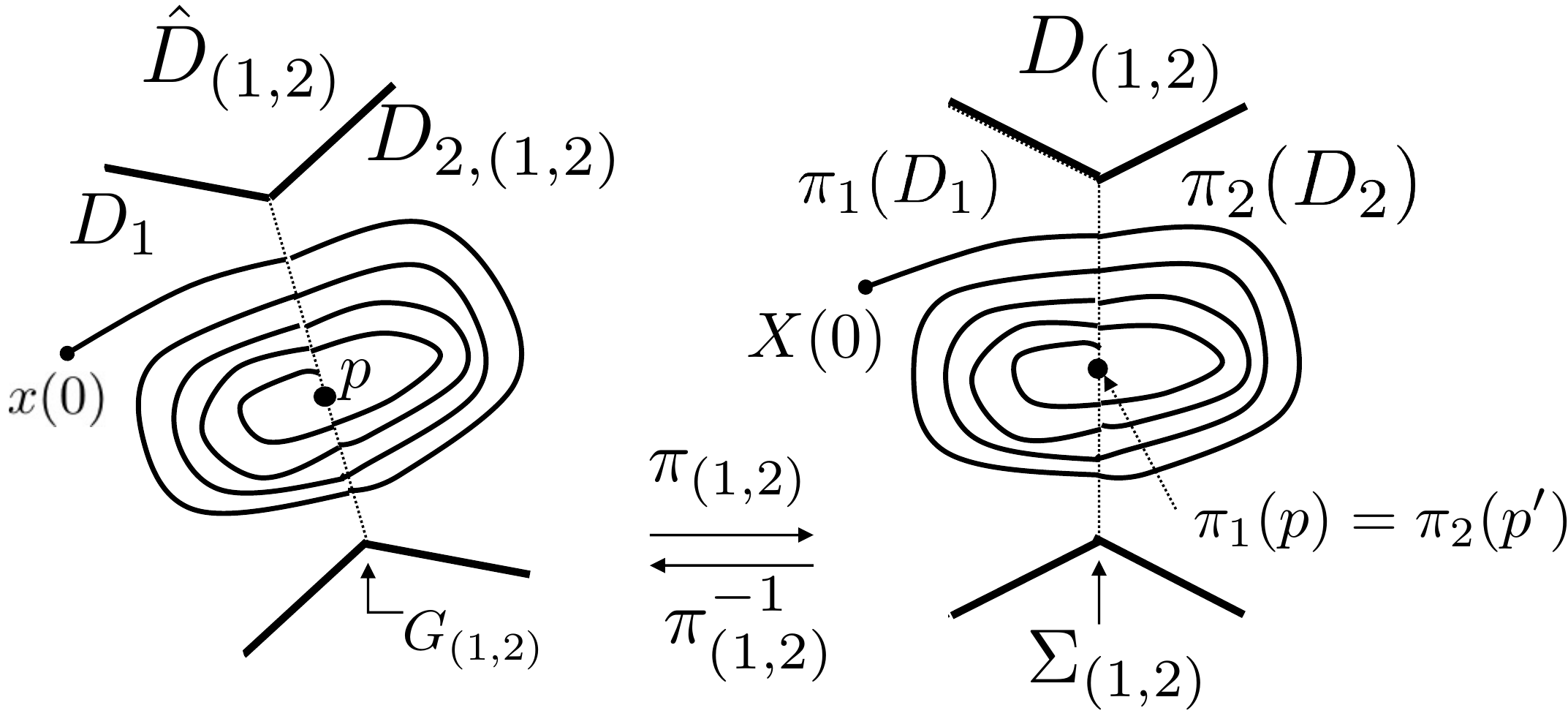}%
  \caption{ The hybrid system from Figure \ref{fig:classic_traj} with the edge $(2,1)$ removed. The set $\hat{D}_{(1,2)}$ is depicted on the left, and $D_{(1,2)} = \M$ is depicted on the right. A hybrid Filippov solution $X$ spirals towards the point $\pi_{1}(p) = \pi_{2}(p')$, crossing $\Sigma_{(1,2)}$ an infinite number of times. The curve $X$ can be constructed by setting $X = \pi_{(1,2)}\circ x$, where $x$ is a Filippov solution of the vector field $f_{(1,2)}$ with initial condition $x(0) = \pi_{(1,2)}^{-1}(X(0))$.}
  \label{fig:hybrid_filippov}
    \vspace{-10pt}
\end{figure}

First, for each $e =(j,j') \in \Gamma$, we define 
\begin{equation}
D_e \colon = int(\pi_j(D_{j})) \bigcup \Sigma_e \bigcup int(\pi_{j'}(D_{j'})),
\end{equation}
which is an open subset of $\M$, and shown in Figure \ref{fig:hybrid_filippov}. In order to construct a coordinate chart for $D_e$, we need to define a set $\hat{D}_e \subset \R^n$ and a bijection from $D_e$ to $\hat{D}_e$. However, we first require a few intermediate constructions. For each $e = (j,j') \in \Gamma$, define $p_e \colon \R^n \to \R^n$ by $p_e(x) = x - \hat{g}_e g_e(x)$, the Euclidian projection onto $\tilde{G}_e$.\footnote{When $\tilde{G}_e$ is a general nonlinear surface, it may not be possible to write down a closed-form expression for the projection onto $\tilde{G}_e$, as our approach requires. } Next, define the map $\bar{R}_e \colon \R^n \to \R^n$ by
\begin{equation}
\bar{R}_e(x) = R_e \circ p_e(x) + \hat{r}_{e} g_e(x),
\end{equation}
and consider the set $D_{j',e}  \colon = \bar{R}_e^{-1}(D_{j'})$. As depicted in Figure \ref{fig:attach}, $D_{j',e}$ is the result of attaching $D_{j'}$ to $D_{j}$, by passing $D_{j'}$ through the map $\bar{R}_e^{-1}$. To understand $\bar{R}_e$, first note that each $x \in \R^n$ may be decomposed as $x = p_e(x) + \hat{g}_eg_e(x)$. Intuitively, $\bar{R}_e$ performs this decomposition, and then takes $p_e(x)$ to $R_e(p_e(x))$ and $\hat{g}_{e} g_e(x)$ to $ \hat{r}_{e} g_e(x)$. Thus, the first term in $\bar{R}_e$ maintains the one-to-one relationship between $\tilde{G}_e$ and $\tilde R_e$ that $R_e|_{\tilde{G}_e}$ defines, while the second term in $\bar{R}_e$ defines a correspondence between the direction that is normal to $\tilde{G}_e$ and the direction that is normal to $\tilde{R}_e$. It is not difficult to verify that $\bar{R}_e$ is a diffeomorphism, and in \cite{arxiv} we give a closed-form representation for its inverse. 

Finally, for each $e =(j,j') \in \Gamma$ we define 
\begin{equation}
\hat{D}_e = int(D_j)\bigcup G_e \bigcup int(D_{j',e})
\end{equation}
and then define the bijection $\pi_e \colon \hat{D}_e \to D_e$ by
\begin{equation}
\pi_e(x) = \begin{cases}
\pi_{j}(x) &\text{if } x \in D_j \\
\pi_{j'}\circ \bar{R}_e(x) & \text{if } x \in D_{j',e}.
\end{cases}
\end{equation}

\begin{theorem}
Let $\hs$ be a hybrid dynamical system. Then we may endow $\M$ with the structure of a smooth topological manifold whose smoothly compatible atlas of coordinate charts is given by $\set{D_e , \pi_e^{-1}}_{e \in \Gamma}$.
\end{theorem}

The details of the proof are given in \cite{arxiv}, but the argument closely follows the proof of \cite[Theorem 9.29]{manifolds2012introduction}. Note that \cite{simic2005towards} and \cite{burden2015model} both showed that it was possible to construct a family of coordinate charts for $\M$. However, in both of these works the coordinate charts are defined implicitly, and thus do not provide explicit representations for portions of $\M$. By providing closed-form representations for our coordinate charts, in Section \ref{sec:hybrid_filippov} we will be able to directly analyze portions of the hybrid Filippov solution using concrete vector fields defined on $\set{\hat{D}_e}_{e \in \Gamma}$.  However, our approach does rely on Assumption \ref{ass:affine}, which is not made in either \cite{simic2005towards} or \cite{burden2015model}. 
\subsection{Relaxed Hybrid Topology}\label{sec:relaxed_topology}

\begin{figure}[t]
\vspace{10pt}
  \centering
  \includegraphics[width= .9\columnwidth, height = .95\textheight, keepaspectratio ]{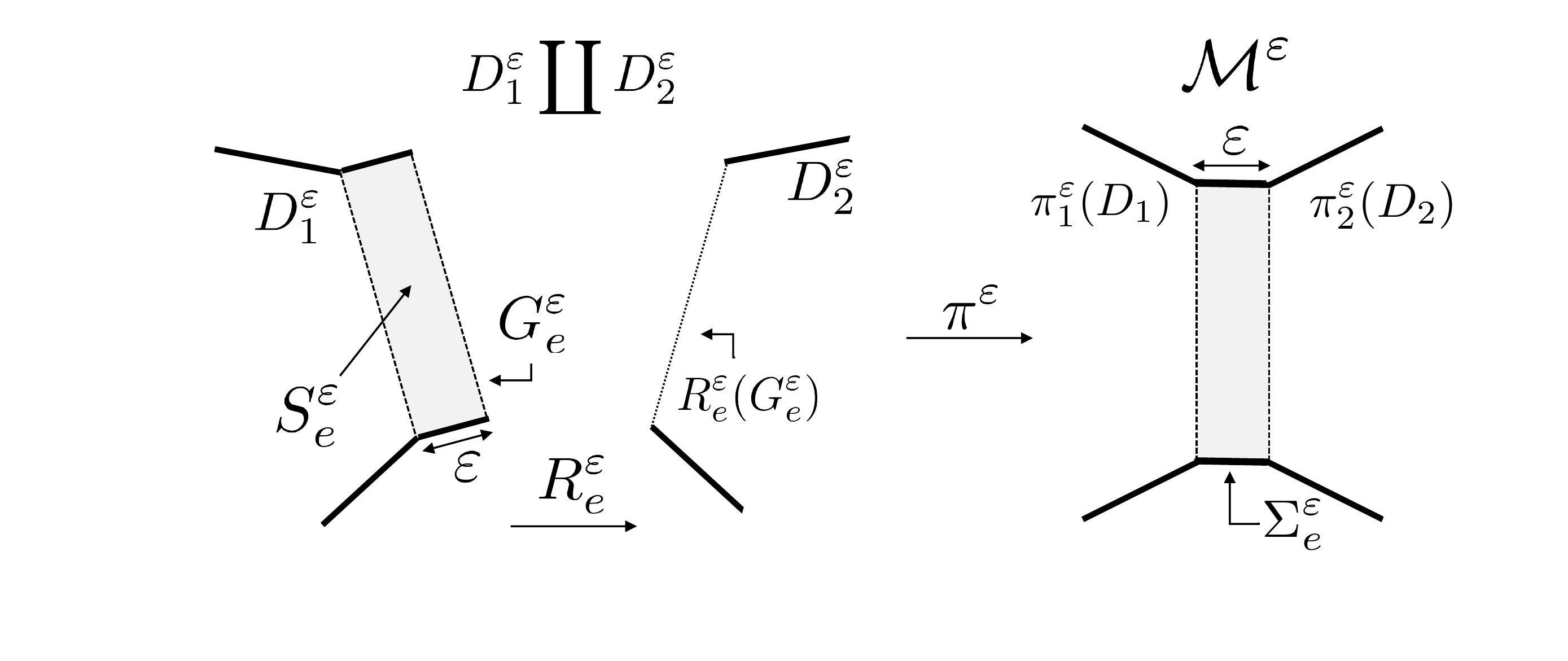}%
  \caption{Construction of the relaxed hybrid quotient space from the disjoint union of the relaxed continuous domains for a bimodal hybrid system with a single edge $e =(1,2)$.}
  \label{fig:relaxed_hybrid_quotient_space}
  \vspace{-10pt}
\end{figure}

We now introduce the relaxed hybrid topology from \cite{burden2015metrization}, which is constructed by attaching an $\epsilon$-thick strip to each of the guard sets of the hybrid system. Tolerances of this sort have been widely used to ensure hybrid models accurately reflect the dynamics of the physical process they are meant to represent \cite{johansson1999regularization}. We also introduce relaxed versions for a number of our previous constructions. This introduction is brief, since many of these objects are quite similar to previous definitions. However, many of these items are depicted in either Figure \ref{fig:relaxed_hybrid_quotient_space} or Figure \ref{fig:relaxed_traj}.

First, for each $e\in \Gamma$ and $\epsilon >0$ we define the \emph{relaxed strip} 
\begin{equation}
S_e^\epsilon := \set{p + \hat{g}_eq \in \R^n \colon p \in G_e \text{ and } q \in \sp{0, \epsilon}}.
\end{equation}
For each $j \in \J$ we then define 
\begin{equation}
{\mc{N}_{j} \colon =\set{e \in \Gamma \colon \exists j' \in \J \text{ s.t. } e =(j,j')}},
\end{equation} the set of edges leaving mode $j$, and then define
the \emph{relaxed domain}  ${D_j^\epsilon =D_j \cup_{e \in \mathcal{N}_j} S_e^\epsilon}$.
Next, for each $e =(j,j') \in \Gamma$ we define the \emph{relaxed guard set} 
\begin{equation}
G_e^\epsilon := \set{x \in S_e^\epsilon  :g_e^\epsilon(x) \colon = \hat{g}_e^Tx - \p{c_e +\epsilon} = 0  },
\end{equation}
and then define the relaxed reset map $R_e^\epsilon \colon \R^n \to \R^n$ by $R_e^\epsilon(x) = R_e(x- \hat{g}_e \epsilon)$. 
Note that $R_e^\epsilon(G_e^\epsilon) = R_e(G_e)$. 

Next, we define $\hat{R}^\epsilon \colon \coprod_{e \in \Gamma} G_e^\epsilon \to \coprod_{j \in \J}D_j^\epsilon$ by $\hat{R}^\epsilon(x) = R_e^\epsilon(x)$ if $x \in G_e^\epsilon$ and then define the \emph{relaxed hybrid quotient space} to be
\begin{equation}
\Me = \frac{\coprod_{j \in \J}D_j^\epsilon}{\Lambda_{\hat{R}^\epsilon}}.
\end{equation}
The construction of the relaxed hybrid quotient space is depicted in Figure \ref{fig:relaxed_hybrid_quotient_space}. For the rest of the paper, we let $d^\epsilon \colon \M^\epsilon \times \M^\epsilon \to \R_{+}$ be a state-space metric on $\M^\epsilon$. It is understood that continuity on $\Me$ is defined with respect to $d^\epsilon$. Letting $\pi^\epsilon \colon \coprod_{j \in \J} D_j^\epsilon \to \Me$ denote the quotient map for $\Me$, for each $j \in \J$ we define the map $\pi_j^\epsilon \colon D_j^\epsilon \to \M^\epsilon$ by $\pi_j^\epsilon(x) = \pi^\epsilon(x \times \set{j})$. Then for each $e =(j,j') \in \Gamma$ we define $\Sigma_e^\epsilon \colon = \pi_{j}^\epsilon(S_e^\epsilon)$. Just as $\Sigma_e$ is a surface separating $\pi_j(D_j)$ from $\pi_{j'}(D_{j'})$ on $\M$, $\Sigma_e^\epsilon$ is an $\epsilon$-thick strip that separates $\pi_j^\epsilon(D_j)$ and $\pi_{j'}^\epsilon (D_{j'})$ on $\Me$. 

Next, we construct a collection of coordinate charts for $\Me$. First, for each $e =(j,j') \in \Gamma$ we define 
\begin{equation}
{D_e^\epsilon \colon = int(\pi_j^\epsilon(D_j)) \bigcup \Sigma_e^\epsilon \bigcup int(\pi_{j'}^\epsilon(D_{j'}))}.
\end{equation}
Subsequently, we define the map ${p_e^\epsilon \colon \R^n \to \R^n}$ by ${p_e^\epsilon = x -\hat{g}_e g_e^\epsilon(x)}$, the Euclidian projection onto the plane containing $G_e^\epsilon$.  We then let ${\bar{R}_e^\epsilon \colon \R^n \to \R^n}$ be such that 
\begin{equation}
{\bar{R}_e^\epsilon(x) = R_e^\epsilon \circ p_e^\epsilon(x) + \hat{r}_e g_e^\epsilon(x)},
\end{equation}
and then define ${D_{j',e}^\epsilon = (\bar{R}_e^{\epsilon})^{-1}(D_{j'})}$. Finally, we define 
\begin{equation}
{\hat{D}_e^\epsilon \colon= int(D_j) \bigcup S_e^\epsilon \bigcup int(D_{j',e}^\epsilon)},
\end{equation}
and let the bijection $\pi_e^\epsilon \colon \hat{D}_e^\epsilon \to D_e^\epsilon$ be defined by
\begin{equation} 
\pi_e^\epsilon(x) = \begin{cases}
\pi_j^\epsilon(x) & \text{if } x \in D_j \cup S_e^\epsilon \\
\pi_{j'}^\epsilon \circ \bar{R}_e^\epsilon(x) & \text{if } x \in D_{j',e}^\epsilon.
\end{cases}
\end{equation}

\begin{theorem}
Let $\hs$ be a hybrid dynamical system. Then for each $\epsilon >0$ we may endow $\M^\epsilon$ with the structure of a smooth topological manifold whose smoothly compatible atlas of coordinate charts is given by $\set{D_e^\epsilon , (\pi_e^\epsilon)^{-1}}_{e \in \Gamma}$.
\end{theorem}

\section{Hybrid Filippov Solutions}\label{sec:hybrid_filippov}
We now apply Filippov's solution concept to our class of hybrid dynamical systems. We begin by defining a single piecewise-smooth vector field on $\M$. 

To construct this vector field, for each $j \in \J$ we define $F_j \colon \pi_j(D_j) \times U \to T\pi_j(D_j)$ by $F_j = D\pi_j \circ f_j$. Intuitively, $F_j$ is simply the representation of the vector field $f_j$ on the set $\pi_j(D_j)$.
\begin{definition}
Let $\hs$ be a hybrid dynamical system. We define the \emph{hybrid vector field} $F \colon \M \times U  \to T\M$ by
\begin{equation}
F(X,u) = F_{j}(X,u) \text{ if } X \in int(\pi_j(D_j)).
\end{equation}
\end{definition}

Note that $F$ is discontinuous and undefined on the surface $\Sigma_{e}$, for each $e \in \Gamma$. Due to Assumption \ref{ass:overlapping_guards}, these surfaces of discontinuity do not intersect. In Section \ref{sec:examples} and \cite{arxiv} we consider several examples where surfaces of discontinuity in the hybrid vector field overlap.

For each $e =(j,j') \in \Gamma$, we define $f_e \colon \hat{D}_e \times U \to \R^n$ by $f_e = D\pi_e^{-1} \circ F$, the coordinate representation of $F$ on the set $\hat{D}_e$. The closed form representation of $f_e$ is given by 
\begin{equation}\label{eq:local_vfield}
f_e(x,u) = \begin{cases}
f_j(x,u) &\text{if } x \in int(D_j) \\
f_{j',e}(x,u) & \text{if } x \in int(D_{j',e}),
\end{cases}
\end{equation}
where $f_{j',e} \colon D_{j',e} \times U \to \R^n$ is defined by
\begin{equation}
f_{j',e}(x,u) = \p{\nabla \bar{R}_e(x)}^{-1} \cdot f_{j'}(\bar{R}_e(x),u).
\end{equation}
Just as $F_{j'}$ is the representation of $f_{j'}$ on the set $\pi_{j'}(D_{j'})$, $f_{j',e}$ is the representation of $f_{j'}$ on the set $D_{j',e}$. The following result is a consequence of Assumptions \ref{ass:lipschitz_vectorfield} and \ref{ass:diff_rmaps} and ensures that the result of Lemma \ref{lemma:existence} applies to $f_e$.
\begin{lemma}\label{lemma:lipschitz1}
For each $e = (j,j') \in \Gamma$ the vector field $f_{j',e}$ is smooth and Lipschitz continuous.
\end{lemma}

We now use the Filippov regularizations for $\set{f_e}_{e \in \Gamma}$ to define the hybrid Filippov regularization for $F$, and subsequently the hybrid Filippov solution. 

\begin{definition}
Let $\hs$ be a hybrid dynamical system. The \emph{hybrid Filippov regularization} of the vector field $F$ is the set-valued map $\widehat{\mc{F}}\sp{F} \colon \M \times U \to \mc{B}\p{T\M}$ where 
\begin{equation}
\widehat{\mc{F}}[F]|_{D_e \times U} = D\pi_e(\mc{F}[f_e]), \ \forall e \in \Gamma.
\end{equation}
For initial condition $X_0 \in \M$ and input $u \in \BVU$, we say that the absolutely continuous curve $X \colon \sp{0,T} \to \M$ is a \emph{hybrid Filippov solution} for this data if $X(0) = X_0$ and
\begin{equation}
 \dot{X}(t) \in \widehat{\mc{F}}\sp{F}(X(t),u(t)) \text{ a.e. } t \in \sp{0,T}.
\end{equation}
\end{definition}

 As shown in Figure \ref{fig:hybrid_filippov}, each portion of a hybrid Filippov solution can be explicitly constructed using a Filippov solution corresponding to one of the piecewise-smooth vector fields $\set{f_e}_{e \in \Gamma}$. Thus, the vector fields $\set{f_e}_{e\in \Gamma}$ can be used to locally assess properties of the hybrid Filippov solution such as stability, controllability, and dependence on initial conditions and parameters \cite{smirnov2002introduction}, \cite{filippov2013differential}. In \cite{arxiv},  we demonstrate how to construct a full hybrid Filippov solution using a boundary value problem similar to the one in Definition \ref{def:classical_hybrid_traj}, in a manner resembling how integral curves are usually constructed on topological manifolds \cite[Chapter 9]{manifolds2012introduction}.
 
Note that hybrid Filippov solutions can cross back and forth across each surface of discontinuity in $F$. Consequently, for each $e =(j,j') \in \Gamma$, it is as if we have implicitly added an extra edge $\bar{e} =(j',j)$ to the hybrid system, where $G_{\bar{e}} = R_e(G_e)$ and $R_{\bar{e}} = R_e^{-1}$. Thus, as depicted in Figure \ref{fig:hybrid_filippov}, we can safely remove one of the edges of the hybrid system depicted in Figure \ref{fig:classic_traj}, and still have the hybrid Filippov solution faithfully capture the dynamics of the system. Note that the hybrid Filippov solution depicted in Figure \ref{fig:hybrid_filippov} can be constructed using a single Filippov solution of the vector field $f_{(1,2)}$, whereas constructing a hybrid execution for this data would have required an infinite number of reset map evaluations.

Finally, we discuss the existence and uniqueness of the hybrid Filippov solution for $\hs$. By Lemmas \ref{lemma:existence} and \ref{lemma:lipschitz1}, for each $e \in \Gamma$, Filippov solutions for the vector field $f_e$ exist on bounded time intervals, up until the solutions leave $\hat{D}_e$. Consequently, hybrid Filippov solutions for $\hs$ exist on bounded time intervals, up until each solution leaves $\M$. The following assumption is analogous to Assumption \ref{ass:transverse}.

\begin{assumption}\label{ass:transverse_hybrid}
Let $\hs$ be a hybrid dynamical system. Then for each $e =(j,j') \in \Gamma$ and each $(x,u) \in G_e \times U$ either $\hat{g}_e^{T} \cdot f_{j}(x,u) >0$ or $\hat{r}_e^T \cdot f_{j'}(R_e(x),u) <0$. 
\end{assumption}

Indeed, by carefully inspecting the terms in \eqref{eq:local_vfield}, one can see that if Assumption \ref{ass:transverse_hybrid} holds $\hs$, then Assumption \ref{ass:transverse} holds for $f_e$, for each $e \in \Gamma$.

\begin{theorem}\label{thm:unique_filippov}
Let $\hs$ be a hybrid dynamical system satisfying Assumption \ref{ass:transverse_hybrid}. Suppose there exists a hybrid Filippov solution $X \colon \sp{0,T} \to \M$ for the data $(X_0,u)  \in  \M \times \BVU$. Then $X$ is the unique hybrid Filippov solution corresponding to this data. 
\end{theorem}

In general, if one is able to verify the uniqueness of the Filippov solutions for the vector fields $\set{f_e}_{e \in \Gamma}$ using a know result, then the uniqueness of the hybrid Filippov solution can be immediately verified.

\section{Relaxed Hybrid Vector Fields}\label{sec:relaxed_hybrid_sys}

Using our previous constructions, we can easily extend the smoothing technique introduced in Section \ref{sec:smooth_systems} to the hybrid vector field $F$. For each $\epsilon >0$, this process will result in a smooth vector field $F^\epsilon \colon \Me \times U \to T\Me$.

First, for each $e =(j,j') \in \Gamma$ and $\epsilon >0$ we define $f_{j',e}^\epsilon \colon D_{j',e}^\epsilon \times U \to \R^n$ by
\begin{equation}
f_{j',e}^\epsilon(x,u) = \p{\nabla \bar{R}_e^{\epsilon}(x)}^{-1} \cdot f_{j'}(\bar{R}_e^\epsilon(x),u),
\end{equation}
which is the representation of $f_{j'}$ on $D_{j',e}^\epsilon$. Note that $f_{j',e}^\epsilon$ is just a translated version of $f_{j',e}$.
\begin{lemma}
For each $e =(j,j') \in \Gamma$ the vector field $f_{j',e}^\epsilon$ is smooth and Lipschitz continuous.
\end{lemma}

We now modify the class of transition functions introduced in Section \ref{sec:smooth_systems}.
\begin{definition}
We say that $\phi \in C^{\infty}(\R, \sp{0,1})$ is a \emph{hybrid transition function} if $i)$ $\phi(a) = 0$ if $a \leq 0$, $ii)$ $\phi(a)= 1$ if $a \geq 1$, and $iii)$ $\phi$ is monotonically increasing on $(0,1)$.
\end{definition}

For the rest of the section, assume a single hybrid transition function $\phi$ has been chosen. Next, for each $e =(j,j') \in \Gamma$ and $\epsilon >0$ we define ${\phi_e^\epsilon \colon \R^n \to \R}$ by ${\phi_e^\epsilon(x) = \phi(\frac{g_e(x)}{\epsilon})}$, and then define the vector field $f_e^\epsilon \colon \hat{D}_e^\epsilon \times U \to \R^n$ by
\begin{equation}
f_e^\epsilon(x,u) = \p{1- \phi_e^\epsilon(x)}f_j(x,u)+  \phi_e^\epsilon(x) f_{j',e}^\epsilon(x,u).
\end{equation}

Note $f_e^\epsilon(x,u) = f_j(x,u)$ if $x \in D_j$, $f_e^\epsilon(x,u) = f_{j',e}^\epsilon(x,u)$ if $x \in D_{j',e}^\epsilon$, and that $f_e^\epsilon(x,u)$ produces a convex combination of $f_j(x,u)$ and $f_{j',e}^\epsilon(x,u)$ when $x \in S_e^\epsilon$. 

\begin{lemma}
Let $\hs$ be a hybrid dynamical system. Then for each $e \in \Gamma$ and $\epsilon>0$ the vector field $f_e^\epsilon$ is smooth.
\end{lemma}

For each $e =(j,j') \in \Gamma$ we then define ${F_e^\epsilon \colon D_e^\epsilon \times U \to TD_e^\epsilon}$ by ${F_e^\epsilon = D\pi_e^\epsilon \circ f_e^\epsilon}$, which smoothly transitions between the dynamics of mode $j$ and the dynamics of mode $j'$ along $\Sigma_e^\epsilon$.

\begin{definition}
Let $\hs$ be a hybrid dynamical system. Then for each $\epsilon>0$ we define the \emph{relaxed hybrid vector field} $F^\epsilon \colon \M^\epsilon \times U \to T\M^\epsilon$ by 
\begin{equation}
F^\epsilon(X,u) = F_e^\epsilon(X,u) \text{ if } X \in D_e^\epsilon.
\end{equation}
For initial condition $X_0 \in \M$ and input $u \in \BVU$ we say that the absolutely continuous curve $X^\epsilon \colon \sp{0,T} \to \Me$ is a \emph{relaxed hybrid trajectory} corresponding to this data if $X^{\epsilon}(0)= X_0$ and
\begin{equation}
\dot{X}^\epsilon(t) = F^\epsilon(X^\epsilon(t),u(t)) \text{ a.e. } \forall t \in \sp{0,T}.
\end{equation} 
\end{definition}

Note that we do not allow relaxed hybrid trajectories to begin on the strips $\set{\Sigma_e^\epsilon}_{e \in \Gamma}$, since they are virtual objects not present in the original hybrid system, instead requiring their initial conditions to be on $\M$, which we regard as a subset of $\M^\epsilon$. The next result follows from the fact that $F^\epsilon$ is represented locally by the smooth vector fields $\set{f_e^\epsilon}_{e \in \Gamma}$. 
\begin{theorem}\label{thm:smooth_vector_field}
Let $\hs$ be a hybrid dynamical system. Then for each $\epsilon > 0$ the vector field $F^\epsilon$ is smooth.
\end{theorem}

\begin{figure}[t]
\vspace{10pt}
  \centering
  \includegraphics[width= .9\columnwidth, height = .95\textheight, keepaspectratio ]{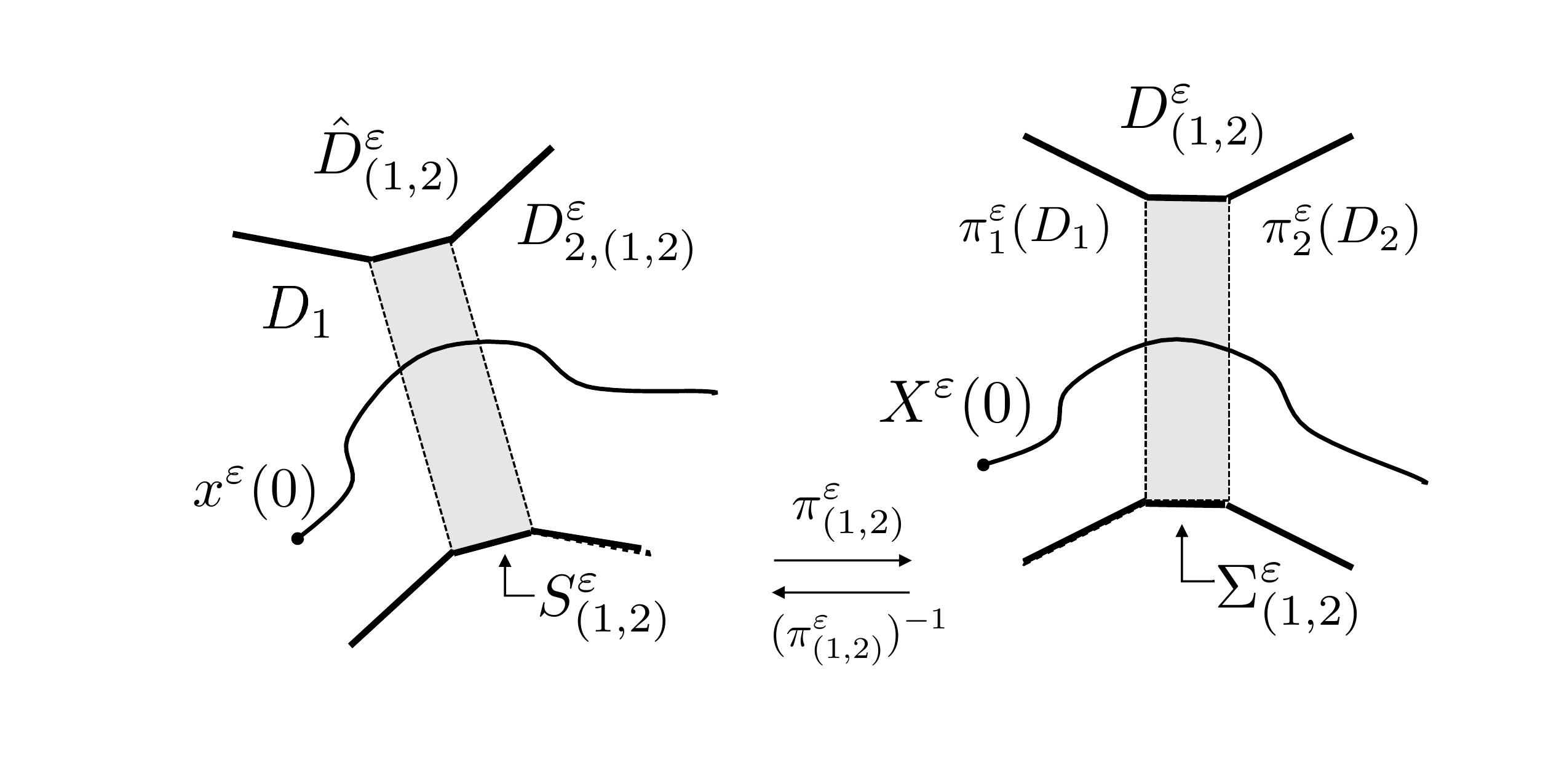}%
  \caption{  A  relaxed hybrid trajectory $X^\epsilon$ with initial condition $X^\epsilon(0)$ flowing on $D_{(1,2)}^\epsilon \subset \M^\epsilon$. This flow can be constructed by setting $X^\epsilon = \pi_{(1,2)}^\epsilon \circ  x^\epsilon$ where $x^\epsilon$ is an integral curve of $f_{(1,2)}^\epsilon$ with initial condition $x^\epsilon(0) = (\pi_{(1,2)}^\epsilon)^{-1}(X^\epsilon(0))$. }
  \vspace{-15pt}
  \label{fig:relaxed_traj}
\end{figure}

It is a fundamental result that the flows of smooth vector fields depend smoothly on inputs and initial conditions \cite[Chapter 4.3]{schattler2012geometric}. The construction of the vector field $F^\epsilon$ can be thought of as a generalization of the \emph{regularization in space} introduced in \cite{johansson1999regularization}, or an approximation to the smoothing techniques discussed in \cite{burden2015model} and \cite{simic2005towards}. As depicted in Figure \ref{fig:relaxed_traj}, each portion of a relaxed hybrid trajectory can be constructed using one of the vector fields $\set{f_e^\epsilon}_{e \in \Gamma}$. In \cite{arxiv}, we demonstrate how to construct a full relaxed hybrid trajectory using a boundary value problem involving the vector fields $\set{f_e^\epsilon}_{e \in \Gamma}$, which can be discretized using standard techniques for integrating stiff differential equations. Finally, we study the limiting behavior of our relaxed solution concept.

\begin{theorem}\label{thm:convergence3}
Let Assumption \ref{ass:transverse_hybrid} hold for hybrid dynamical system $\hs$. Let $(X_0, u) \in \M \times \BVU$, and let $X \colon \sp{0,T} \to \M$ be a hybrid Filippov solution corresponding to this data, guaranteed to be unique by Theorem \ref{thm:unique_filippov}. For each $\epsilon>0$, let $X^\epsilon \colon \sp{0,T} \to \M^\epsilon$ be the corresponding relaxed hybrid trajectory. Then $\exists C > 0$ and $\epsilon_0 >0$ such that for each $\epsilon \leq \epsilon_0$ and $t \in \sp{0,T}$ 
\begin{equation}
d^\epsilon\p{X(t), X^\epsilon(t)} \leq C \epsilon.
\end{equation}
\end{theorem}

Here, as in \cite{burden2015metrization}, we have abused notation and regarded a hybrid Filippov solution $X$ as a curve on $\Me$. Next, we demonstrate that our relaxations always converge uniformly to a unique, well-defined limit even in cases where the hybrid Filippov solution may be non-unique.
\begin{theorem}\label{thm:convergence4}
Let $\hs$ be a hybrid dynamical system. Fix $(X_0, u) \in \M \times \BVU$,  and for each $\epsilon >0$ let $X^\epsilon \colon \sp{0,T} \to \Me$ be the corresponding relaxed hybrid trajectory. Then there exists an absolutely continuous curve $X^0 \colon \sp{0,T} \to \M$ such that for each $t \in \sp{0,T}$
\begin{equation}
\lim_{\epsilon \to 0} d^\epsilon \p{X^0(t), X^\epsilon(t)} = 0.
\end{equation}
\end{theorem}

Due to the uniqueness of this limit, we have found it convenient to think of hybrid dynamics as the limit of our relaxations. Further work is required to determine whether the limit in Theorem \ref{thm:convergence4} is in fact a hybrid Filippov solution.

\section{Modeling Example}\label{sec:examples}
We breifly examine an example where surfaces of discontinuity in the hybrid vector field overlap. Further details are provided in \cite{arxiv}, where additional examples are considered.

Consider a ball which is bouncing vertically and loses a fraction of its energy during each impact. We model the ball with two continuous states $x =(x_1,x_2)^T$, where $x_1$ is the height of the ball above the ground and $x_2$ is the velocity of the ball. When airborne, the continuous states evolve according to $\deriv{}{t}(x_1,x_2)^T =  [x_2, -g]^T$, where $g >0$ is the gravitational constant. When the ball hits the ground, the velocity is reset according to $x_2 \to - c x_2$, where $c \in(0,1]$ is the \emph{coefficient of restitution}. It is well known \cite{johansson1999regularization} that for $c \in (0,1)$ the ball bounces an infinite number of times by some finite time $t_{\infty}$, thus classical constructions of hybrid executions for the bouncing ball are necessarily Zeno.

We model the bouncing ball with four discrete modes: $\J_{bb} = \set{1,2,3,4}$. For $j \in \set{1,3}$ we define $D_{j} = \set{(x_1,x_2) \in \R^2 \colon x_1 \geq 0, x_2 \geq 0}$, and for $j \in \set{2,4}$ we define $D_{j} =  \set{(x_1,x_2) \in \R^2 \colon x_1 \geq 0, x_2 \leq 0}$. Thus in modes 1 and 3 the ball is moving upwards, and in modes $2$ and $4$ the ball is moving downwards. Our model has four edges: $\Gamma_{bb} = \set{(1,2), (2,3), (3,4), (4,1)}$. For $e \in \set{(1,2),(3,4)}$ we define $G_e =\set{(x_1,x_2)\in \R^2 \colon x_1 \geq 0, x_2  =0}$ and $R_e(x) =x$. For $e \in\set{(2,3),(4,1)}$ we define $G_e =\set{(x_1,x_2) \in \R^2 \colon x_1 =0 , x_2  \leq 0}$ and $R_e(x) =[x_1, -c x_2]^{T}$. Thus edges $(1,2)$ and $(3,4)$ are triggered when the ball is at the apex of its flight, and edges $(2,3)$ and $(3,4)$ are triggered when the ball hits the ground.

\begin{figure}[t]
\vspace{10pt}
  \centering
  \includegraphics[width= .9\columnwidth, keepaspectratio ]{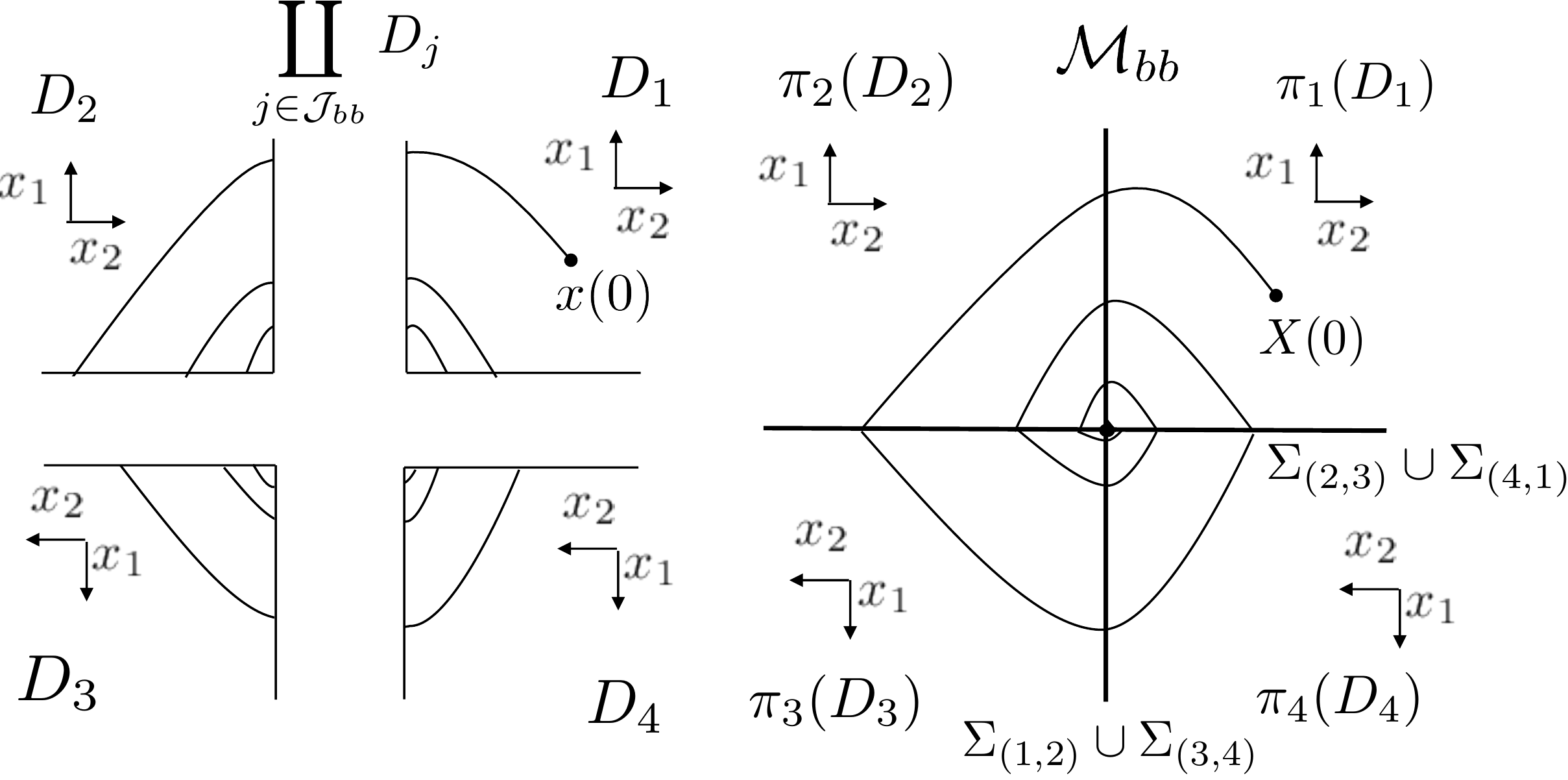}%
  \caption{ A hybrid execution $x$ for the bouncing ball evolves on the disjoint union of the continuous domains, and the corresponding hybrid Filippov solution $X$ evolves on $\M_{bb}$. In each domain, the axis denotes the orientation of the two continuous states.}
  \vspace{-10pt}
  \label{fig:bb}
\end{figure}

The hybrid quotient space for the bouncing ball, $\M_{bb}$, is depicted in Figure \ref{fig:bb}. In \cite{arxiv}, we demonstrate how to construct a single chart for $\M_{bb}$ wherein each continuous domain is sent to one of the quadrants of $\R^2$. The hybrid vector field has two overlapping surfaces of discontinuity: $\Sigma_{(1,2)} \cup \Sigma_{(3,4)}$ and $\Sigma_{(2,3)} \cup \Sigma_{(4,1)}$. A hybrid Filippov solution for the bouncing ball is depicted in Figure \ref{fig:bb}. The trajectory crosses the surfaces of discontinuity an infinite number of times before coming to rest at the origin for each $t \geq t_{\infty}$, the hybrid Filippov solution naturally extending the trajectory past the Zeno point. As addressed in \cite{arxiv}, we extend the smoothing technique from \cite{da2018slow} to relax the dynamics of the bouncing ball.

\section{Conclusion and Future Work}
We employed Filippov's solution concept for differential equations with discontinuous right-hand sides to describe the trajectories of a class of hybrid systems which display discrete jumps in the continuous state. We then introduced a family of smooth vector fields which can be used to approximate these dynamics in the numerical setting, using existing methods. Further work is required to characterize these two new solution concepts, and to control their trajectories using existing analytic and computational tools.

\bibliographystyle{IEEEtran}
\bibliography{paper}
\end{document}